# NON-EXISTENCE OF CERTAIN SEMISTABLE ABELIAN VARIETIES

ARMAND BRUMER AND KENNETH KRAMER

## 1. INTRODUCTION

If $A$ is an abelian variety defined over $\mathbb{Q}$, its Tate module $\mathbb{T}_\ell(A)$ affords an $\ell$-adic representation $\rho$ of the absolute Galois group $G_\mathbb{Q}$. But under suitable constraints on ramification, such a representation cannot exist. Thus, Fontaine [Fo] proved there do not exist abelian varieties over $\mathbb{Z}$, *i.e.* with everywhere good reduction. Fontaine speculated [Fo, Rem. 3.4.7] that the same methods might rule out semistable varieties defined over $\mathbb{Q}$ with good reduction outside one small prime. Fontaine's methods were extended by Joshi [Jo] to prove non-existence of certain crystalline mod-$\ell$ representations $\bar{\rho}: G_\mathbb{Q} \to \mathrm{GL}_n(\mathbb{F}_\ell)$. Recently, Schoof [Sc] announced a further extension of these ideas to the study of abelian varieties with everywhere good reduction over larger number fields. In this note, we verify that there do not exist semistable abelian varieties over $\mathbb{Q}$ with good reduction outside one prime $p$, for $p \in \{2, 3, 5, 7\}$. For these primes, it is well-known that there are elliptic curves with good reduction outside $p$, but of course not semistable, and there are semistable elliptic curves with good reduction outside $p = 11$. See also [MS] for abelian varieties of dimension 2 with good reduction outside 2, not semistable.

We may summarize our version of these methods as follows. Restrictions on the ramification in $\rho$ imply an upper bound for the discriminant of the $\ell$-division field $L = \mathbb{Q}(A[\ell])$. Odlyzko's work [Od] on discriminants then limits the possible fields $L$. When we have sufficiently good control over $L$, we can construct a chain of isogenies

$$A \to A' \to A'' \to \ldots$$

involving arbitrarily many non-isomorphic abelian varieties. This contradicts the Shafarevich conjecture, as proved by Faltings (see [Fa] and [Wu, Thm. 3.1 ff]).

While Fontaine and Joshi are mainly concerned with the restriction of $\rho$ to a decomposition group over $\ell$, we make more careful use of the local behavior of $A$ over the completion $\mathbb{Q}_p$ at the bad prime $p$. In particular, we introduce an invariant (see §2) that we call the effective stage of $p$-adic inertia acting on $\mathbb{T}_\ell(A)$. This invariant is related to the group of connected components in the special fiber of the Néron model for $A$ at $p$, and is used to guarantee that the varieties in our chain of isogenies are not isomorphic. Furthermore, our approach sometimes applies when $\ell = 2$. See

Proposition 4.3 for an application to the non-existence of abelian varieties with small 2-division fields.

It is convenient to establish some notation here. In general, $A$ denotes an abelian variety defined over a field $K$ of characteristic 0, and $G_K = \text{Gal}(\bar{K}/K)$ is the Galois group of a fixed algebraic closure over $K$. We write $\rho_K$ for the $\ell$-adic representation of $G_K$ afforded by $\mathbb{T}_\ell(A)$. Let $\hat{A}$ denote the dual abelian variety of $A$. The Weil pairing induces a Galois-equivariant perfect pairing

$$(1.1) \qquad e_n : A[\ell^n] \times \hat{A}[\ell^n] \to \boldsymbol{\mu}_{\ell^n},$$

as a consequence of which $\boldsymbol{\mu}_{\ell^n} \subset K(A[\ell^n])$. If $\lambda : A \to \hat{A}$ is a polarization and $P, Q \in A[\ell^n]$, we put $e_n^\lambda(P, Q) = e_n(P, \lambda(Q))$. Passing to the limit gives rise to a perfect pairing $e_\infty : \mathbb{T}_\ell(A) \times \mathbb{T}_\ell(\hat{A}) \to \mathbb{Z}_\ell(1)$, and a pairing $e_\infty^\lambda : \mathbb{T}_\ell(A) \times \mathbb{T}_\ell(A) \to \mathbb{Z}_\ell(1)$. Recall that $\lambda$ induces an injection $\mathbb{T}_\ell(\lambda) : \mathbb{T}_\ell(A) \to \mathbb{T}_\ell(\hat{A})$ whose cokernel is finite in general, and trivial if and only if the degree of $\lambda$ is prime to $\ell$.

## 2. Local considerations

In this section, we suppose that $K$ is a non-archimedean local field of characteristic 0, with valuation $v_K$ and perfect residue field $k$ of characteristic $p$. We denote by $A$ an abelian variety of dimension $d$ defined over $K$, with semistable bad reduction. To fix notation, recall Grothendieck's decomposition [Gr, §2.5] of $\mathbb{T}_\ell(A)$, assuming $\ell$ is a prime different from $p$. The connected component of the identity $\mathcal{A}_k^0$ of the special fiber $\mathcal{A}_k$ of the Néron model of $A$ admits a decomposition

$$0 \to \mathcal{T} \to \mathcal{A}_k^0 \to \mathcal{B} \to 0$$

in which $\mathcal{T}$ is a torus and $\mathcal{B}$ an abelian variety defined over $k$. Let $\dim \mathcal{T} = t$ and $\dim \mathcal{B} = a$, with $t + a = d$. Write $\Phi_A = \mathcal{A}_k/\mathcal{A}_k^0$ for the group of connected components.

Let $K^{nr}$ be the maximal unramified extension of $K$ inside $\bar{K}$, and denote the inertia group by $\mathcal{I} = \text{Gal}(\bar{K}/K^{nr})$. Put $\mathcal{M}_1 = \mathcal{M}_1(A)$ for the submodule of $\mathbb{T}_\ell(A)$ fixed by $\mathcal{I}$, and let $\mathcal{M}_2$ be the subspace of $\mathbb{T}_\ell(A)$ orthogonal to $\mathcal{M}_1(\hat{A})$ under the $e_\infty$-pairing. According to the Igusa-Grothendieck theorem, semistability of $A$ is equivalent to the containment $\mathcal{M}_2 \subset \mathcal{M}_1$. In that case, we have the decomposition

$$(2.1) \qquad \mathbb{T}_\ell(A) \supset \mathcal{M}_1 \supset \mathcal{M}_2 \supset 0,$$

in which the successive quotients are torsion-free $\mathbb{Z}_\ell$-modules. We may identity $\mathcal{M}_2 \simeq \mathbb{T}_\ell(\mathcal{T})$ and $\mathcal{M}_1/\mathcal{M}_2 \simeq \mathbb{T}_\ell(\mathcal{B})$.

It is well-known [ST, appendix] that $\mathcal{I}$ acts on $\mathbb{T}_\ell(A)$ through its maximal pro-$\ell$ quotient. Moreover, if $g \in \mathcal{I}$, then $(g-1)(\mathbb{T}_\ell(A)) \subset \mathcal{M}_2$. Indeed, using the fact that $g$ acts trivially on $\mathbb{Z}_\ell(1)$ and $\mathcal{M}_1(\hat{A})$, we have

$$e_\infty(g(x) - x, y) = e_\infty(x, y)^{g-1} = 1$$

for all $x \in \mathbb{T}_\ell(A)$ and $y \in \mathcal{M}_1(\hat{A})$. Hence $g(x) - x \in \mathcal{M}_2$. In fact, we have a homomorphism $\mathcal{I} \to \mathrm{Hom}(\mathbb{T}_\ell(A)/\mathcal{M}_1, \mathcal{M}_2)$ induced by $g \mapsto N_g = g - 1$. With respect to the decomposition $\mathbb{T}_\ell(A) \simeq \mathcal{M}_2 \oplus \mathcal{M}_1/\mathcal{M}_2 \oplus \mathbb{T}_\ell(A)/\mathcal{M}_1$ derived from (2.1), we may represent $\rho_K(g)$ in the form

$$\begin{pmatrix} 1_t & 0 & N_g \\ 0 & 1_{2a} & 0 \\ 0 & 0 & 1_t \end{pmatrix},$$

where $1_n$ denotes the $n \times n$ identity. The criterion of Néron-Ogg-Shafarevich guarantees that there is a minimal integer $n \geq 1$ such that the restriction of $\mathcal{I}$ to the $\ell^n$-division field $K(A[\ell^n])$ acts non-trivially.

**Definition.** Assume $A$ has semistable bad reduction, and let $i(A, \ell, v_K)$ denote the minimal integer $n \geq 1$ such that the restriction of $\mathcal{I}$ to $K(A[\ell^n])$ is not trivial. We call $i(A, \ell, v_K)$ the effective stage of inertia acting on $\mathbb{T}_\ell(A)$.

Suppose $i(A, \ell, v_K) = n_0$, and fix a topological generator $\sigma$ for the maximal pro-$\ell$ quotient of $\mathcal{I}$. An equivalent formulation is that $n_0$ is the minimal integer such that $N_\sigma \not\equiv 0 \pmod{\ell^{n_0}}$. This does not depend on the choice of $\sigma$. For $n \geq n_0$, it is clear that $K^{nr}(A[\ell^n])$ is the unique tamely ramified extension of $K^{nr}$ of degree $\ell^{n-n_0+1}$. Hence the different ideal $\mathfrak{D}$ of $K(A[\ell^n])/K$ satisfies

(2.2) $$v_K(\mathfrak{D}) = \ell^{n-n_0+1} - 1.$$

Let $\mathcal{N} = N_\sigma(\mathbb{T}_\ell(A)/\mathcal{M}_1)$ be the image of $N_\sigma$. In fact $\mathcal{N}$ does not depend on the choice of $\sigma$. It is known that $\det(N_\sigma) \neq 0$, or equivalently, $\mathcal{N}$ has finite index in $\mathcal{M}_2$. Indeed, according to [Ed, Remark 2.6], the $\ell$-Sylow subgroup of $\Phi_{\hat{A}}(\bar{k})$ is isomorphic as a Galois module for $\mathrm{Gal}(K^{nr}/K) \simeq G_k$ to the Tate twist $(\mathcal{M}_2/\mathcal{N})(-1)$. Therefore

(2.3) $$\mathrm{ord}_\ell(\Phi_{\hat{A}}(\bar{k})) = \mathrm{ord}_\ell(\det(N_\sigma)).$$

To examine variations in the effective stage of inertia under isogeny, consider a $K$-isogeny $\varphi : A \to A'$. Write $\mathcal{M}'_2 = \mathcal{M}_2(A')$, $\mathcal{M}'_1 = \mathcal{M}_1(A')$, and $N'_\sigma : \mathbb{T}_\ell(A')/\mathcal{M}'_1 \to \mathcal{M}'_2$ for the corresponding constructions arising from $A'$. Denote by $\bar{\mathcal{M}}_1$ and $\bar{\mathcal{M}}_2$ the projections of $\mathcal{M}_1$ and $\mathcal{M}_2$ to $A[\ell]$.

**Lemma 2.4.** *Suppose $\kappa$ is a proper $G_K$-submodule of $A[\ell]$ and let $\varphi : A \to A'$ be the $K$-isogeny with kernel $\kappa$. Then*

$$\mathrm{ord}_\ell(\Phi_{\hat{A}'}(\bar{k})) + \dim \kappa/(\kappa \cap \bar{\mathcal{M}}_1) = \mathrm{ord}_\ell(\Phi_{\hat{A}}(\bar{k})) + \dim \kappa \cap \bar{\mathcal{M}}_2.$$

*If $\bar{\mathcal{M}}_2 \subset \kappa$, then $i(A', \ell, v_K) \geq i(A, \ell, v_K)$. If $\bar{\mathcal{M}}_2 \subset \kappa \subset \bar{\mathcal{M}}_1$, then $i(A', \ell, v_K) = i(A, \ell, v_K) + 1$.*

**Proof**. We are content to outline the argument. Let $\varphi'$ be the isogeny from $A'$ whose kernel is $\varphi(A[\ell])$, so that $\varphi' \circ \varphi = [\ell]_A$ is multiplication by $\ell$ on $A$ and $\varphi \circ \varphi' = [\ell]_{A'}$. We have the containments $\mathbb{T}_\ell(\varphi)(M_1) \subset M_1'$ and $\mathbb{T}_\ell(\varphi')(M_2') \subset M_2$. The isogenies $\varphi$ and $\varphi'$ induce maps $\varphi_*$ and $\varphi'_*$ making the following diagram commutative:

(2.5)
$$\begin{array}{ccc} \mathbb{T}_\ell(A)/\mathcal{M}_1 & \xrightarrow{\varphi_*} & \mathbb{T}_\ell(A')/\mathcal{M}_1' \\ \ell N_\sigma \downarrow & & N'_\sigma \downarrow \\ \mathcal{M}_2 & \xleftarrow{\varphi'_*} & \mathcal{M}_2' \end{array}$$

The maps $\varphi_*$ and $\varphi'_*$ are injective. Furthermore, $\dim \operatorname{Coker}(\varphi_*) = \dim \kappa - \dim \kappa \cap \bar{\mathcal{M}}_1$, and $\dim \operatorname{Coker}(\varphi'_*) = \dim \bar{\mathcal{M}}_2 - \dim \kappa \cap \bar{\mathcal{M}}_2$. In view of (2.3), our dimension formula can now be verified by taking the determinant of the relation $\varphi'_* N'_\sigma \varphi_* = \ell N_\sigma$.

Put $n = i(A, \ell, v_K)$, so $N_\sigma \equiv 0 \pmod{\ell^{n-1}}$. If $\bar{\mathcal{M}}_2 \subset \kappa$, we may use the fact that $\varphi_*$ is an isomorphism to show $N'_\sigma \equiv 0 \pmod{\ell^{n-1}}$. Hence $i(A', \ell, v_K) \geq n$. Finally, if $\bar{\mathcal{M}}_2 \subset \kappa \subset \bar{\mathcal{M}}_1$, then both $\varphi_*$ and $\varphi'_*$ are isomorphisms, and the equality $i(A', \ell, v_K) = i(A, \ell, v_K) + 1$ easily follows. ∎

The next two lemmas are not essential for the rest of our argument, but may be of interest if one wishes to work only with principally polarized abelian varieties. Assume for the moment that $A$ admits a principal polarization $\lambda$. Suppose $\kappa \subset A[\ell]$ is maximal isotropic for the perfect pairing $e_1^\lambda : A[\ell] \times A[\ell] \to \boldsymbol{\mu}_\ell$, in the sense that $\kappa = \kappa^\perp$, and let $\varphi : A \to A'$ be the isogeny whose kernel is $\kappa$. Then [Mi, Prop. 16.8] the polarization $\ell\lambda$ induces a principal polarization on $A'$. The following lemmas allow us to construct such maximal isotropic subspaces. Write $E$ (resp. $E_\infty$) for the maximal unramified extension of $K$ inside $L = K(A[\ell])$ (resp. inside $L_\infty = K(A[\ell^\infty])$). Note that $G_K$ acts on $\mathcal{M}_1$, $\mathcal{M}_2$, and $\mathcal{M}_1/\mathcal{M}_2$ via the quotient $\operatorname{Gal}(E_\infty/K)$. Let $\zeta_\infty$ be a generator for $\mathbb{Z}_\ell(1) = \varprojlim \boldsymbol{\mu}_{\ell^n}$, and let $\zeta$ be the projection of $\zeta_\infty$ to $\boldsymbol{\mu}_\ell$.

**Lemma 2.6.** *Assume $A$ admits a principal polarization $\lambda$ defined over $K$. Suppose $\ell$ is odd and there exists an element $\tau$ of order 2 in $\operatorname{Gal}(E/K)$ such that $\tau(\zeta) = \zeta^{-1}$. Then there exists a maximal isotropic subspace $\kappa$ of $A[\ell]$ such that $G_K$ acts on $\kappa$ and $\bar{\mathcal{M}}_2 \subset \kappa \subset \bar{\mathcal{M}}_1$.*

**Proof**. Since $\operatorname{Gal}(L_\infty/L)$ is a pro-$\ell$ group, so is $\operatorname{Gal}(E_\infty/E)$. Hence $\tau$ lifts to an element $\tau_\infty$ of order 2 in $G(E_\infty/K)$ such that $\tau_\infty(\zeta_\infty) = \zeta_\infty^{-1}$. The pairing $e_\infty^\lambda$ induces a perfect pairing
$$e : \mathcal{M}_1/\mathcal{M}_2 \times \mathcal{M}_1/\mathcal{M}_2 \to \mathbb{Z}_\ell(1).$$
Decompose $\tilde{V} = \mathcal{M}_1/\mathcal{M}_2$ into eigenspaces $\tilde{V}^\epsilon = \{v \in \tilde{V} \,|\, \tau_\infty(v) = \epsilon v\}$ for $\epsilon = \pm 1$. Let us verify that $\tilde{V}^\epsilon \subset (\tilde{V}^\epsilon)^\perp$. Indeed, if $v, w \in \tilde{V}^\epsilon$, we have
$$e(v,w)^{-1} = e(v,w)^{\tau_\infty} = e(v^{\tau_\infty}, w^{\tau_\infty}) = e(\epsilon v, \epsilon w) = e(v, w).$$

Hence $e(v, w) = 1$. Then by a standard argument, $\operatorname{rank} \tilde{V}^+ = \operatorname{rank} \tilde{V}^- = \frac{1}{2} \operatorname{rank} \tilde{V} = a$, and $\tilde{V}^\epsilon = (\tilde{V}^\epsilon)^\perp$. Define the pseudo-eigenspaces

$$V^\epsilon = \{m \in \mathcal{M}_1 \,|\, \text{the coset } m + \mathcal{M}_2 \text{ is in } \tilde{V}^\epsilon\}.$$

It easily follows that $V^\epsilon \subset (V^\epsilon)^\perp$ under the pairing $e_\infty^\lambda$. But $(V^\epsilon)^\perp / V^\epsilon$ is torsion-free, and $\operatorname{rank} V^\epsilon = \operatorname{rank} \tilde{V}^\epsilon + \operatorname{rank} \mathcal{M}_2 = a + t = \frac{1}{2} \operatorname{rank} \mathbb{T}_\ell(A)$. Hence $V^\epsilon = (V^\epsilon)^\perp$. We may choose $\kappa$ to be the projection of $V^+$ or $V^-$ to $A[\ell]$. ∎

**Lemma 2.7.** *Assume $A$ admits a principal polarization $\lambda$ defined over $K$. Suppose $\ell = 2$ and $\operatorname{Gal}(E/K)$ is a 2-group. Then there exists a maximal isotropic subspace $\kappa$ of $A[\ell]$ such that $G_K$ acts on $\kappa$ and $\bar{\mathcal{M}}_2 \subset \kappa \subset \bar{\mathcal{M}}_1$.*

**Proof.** Suppose, quite generally, that $q$ is a power of a prime $\ell$ and $V$ is a vector space of dimension $2n$ over $\mathbb{F}_q$ admitting a perfect symplectic pairing. Then the order of the set $S$ of maximal isotropic subspaces of $V$ is $\prod_{i=1}^n (1 + q^i)$. Suppose $H$ is an $\ell$-group and the pairing is $H$-equivariant. Then $H$ acts on $S$ with a fixed point; i.e. there exists a maximal isotropic subspace $W$ of $V$ such that $W$ is a module for $H$.

For the present lemma, consider $V = \bar{\mathcal{M}}_1 / \bar{\mathcal{M}}_2$, upon which $e_1^\lambda$ induces a perfect symplectic pairing. Note that $G_K$ acts on $\bar{\mathcal{M}}_1$, $\bar{\mathcal{M}}_2$ and $V$ through the 2-group $H = \operatorname{Gal}(E/K)$. Let $W$ be a maximal isotropic subspace of $V$ which is an $H$-module as constructed above. Take $\kappa = \{m \in \bar{\mathcal{M}}_1 \,|\, \text{the coset } m + \bar{\mathcal{M}}_2 \text{ is in } W\}$. ∎

In view of various competing notations, we take this opportunity to standardize our ramification numbering, following [Se, Ch. IV] rather than [Fo]. In general, if $L/K$ is a Galois extension of local fields with Galois group $G$, and $\pi_L$ is a prime element of $L$, the (lower) ramification groups are defined by

$$G_n = \{\sigma \in G \,|\, v_L(\sigma(\pi_L) - \pi_L) \geq n + 1\}.$$

For $m \leq u \leq m+1$, the Herbrand function to the upper numbering is given by

$$\varphi_{L/K}(u) = \frac{1}{g_0}(g_1 + \cdots + g_m + (u - m)g_{m+1}),$$

where $g_n = |G_n|$. By definition $G^{n'} = G_n$, where $n' = \varphi_{L/K}(n)$. In terms of this upper numbering, we may restate [Fo, Thm. A] as follows. Let $B[p^n]$ be a finite flat commutative group scheme over the ring of integers $\mathcal{O}_K$ annihilated by $p^n$ for $n \geq 1$. In particular, $B[p^n]$ could be the kernel of multiplication by $p^n$ on an abelian variety $B$ with good reduction. Put $e_K$ for the absolute ramification index of $K$ and let $G = \operatorname{Gal}(K(B[p^n])/K)$.

(2.8) $\quad\quad\quad\quad$ If $u > e_K\left(n + \dfrac{1}{p-1}\right) - 1$ then $G^u = \{1\}$.

Furthermore [Fo, Thm. 1], if the $p$-adic valuation is extended to $L_n = K(B[p^n])$ so that $v_p(p) = 1$, then the different $\mathfrak{D}_{L_n/K}$ satisfies

$$v_p(\mathfrak{D}_{L_n/K}) < n + \frac{1}{p-1}. \tag{2.9}$$

The following lemma will used later, in conjunction with (2.8), to control the conductor of certain abelian extensions. Let $U_K$ be the unit group of $K$ and $U_K^{(n)} = \{x \in U_K \,|\, v_K(x-1) \geq n\}$.

**Lemma 2.10.** *Suppose $L/K$ is Galois extension with $G = \mathrm{Gal}(L/K)$. Consider intermediate fields $L \supset E \supset F \supset K$ such that $E/K$ is Galois and $E/F$ is an abelian $p$-group. Put $e_{tame}$ for the tame ramification degree of $E/K$. Assume $G^u = \{1\}$ for all $u > 1/e_{tame}$. Then the conductor of $E/F$ is at most 2, and the normic subgroup $N_{E/F}(U_E)$ contains $U_F^{(2)}$.*

**Proof.** Let $\bar{G} = \mathrm{Gal}(E/K)$ and $\bar{g}_n = |\bar{G}_n|$. In particular, $\bar{g}_0 = \bar{g}_1 e_{tame}$. Put $u_0 = \varphi_{E/K}(1+\epsilon)$, where $\varphi_{E/K}$ is the Herbrand function and $0 < \epsilon < 1$. Then

$$u_0 = (\bar{g}_1 + \epsilon \bar{g}_2)/\bar{g}_0 > 1/e_{tame},$$

so $G^{u_0} = \{1\}$ by hypothesis. But the upper ramification numbering behaves well under passage to the quotient group $\bar{G}$. Hence $\bar{G}_{1+\epsilon} = \bar{G}^{u_0}$ also is trivial. Let $\bar{H} = \mathrm{Gal}(E/F)$. Since the lower numbering behaves well with respect to subgroups, we have $\bar{H}_{1+\epsilon} = \bar{G}_{1+\epsilon} \cap \bar{H} = \{1\}$. By class field theory [Se, Ch. XV, §2, Cor. 2 to Thm. 1], the conductor of $E/F$ is at most $\varphi_{E/F}(1) + 1 = 2$ and the units in $U_F^{(2)}$ are norms. ∎

As an immediate consequence of the bounds (2.8) and (2.9), we have some control over the group of $p$-power roots of unity contained in the $p^n$-division field of an abelian variety $B$ with good reduction.

**Proposition 2.11.** *Suppose $K/\mathbb{Q}_p$ is unramified and let $B$ be an abelian variety over $K$ with good reduction. Let $L_n = K(B[p^n])$ and $F_n = K(\boldsymbol{\mu}_{p^n})$. If $p$ is odd, $L_n \cap F_\infty = F_n$. If $p = 2$, then $F_n \subset L_n \cap F_\infty \subset F_{n+1}$. Furthermore, if $L_2 \cap F_\infty = F_2$, then $L_n \cap F_\infty = F_n$ for all $n \geq 2$.*

**Proof.** On the one hand, the cyclotomic field $F_n$ is contained in $L_n$. On the other hand, since $K/\mathbb{Q}_p$ is unramified, the different of $F_{n+1}$ satisfies

$$v_p(\mathfrak{D}_{F_{n+1}/K}) = (n+1) - \frac{1}{p-1}.$$

If $p$ is odd, this already exceeds the bound in (2.9) on $\mathfrak{D}_{L_n/K}$, so $F_{n+1} \not\subset L_n$. Similarly, if $p = 2$, then $F_{n+2} \not\subset L_n$

Suppose $p = 2$. To verify that $L_1 \cap F_\infty \subset F_2 = K(\boldsymbol{\mu}_4)$, it suffices to show $\sqrt{\pm 2} \notin L_1$. If $\theta = \sqrt{2}$ or $\sqrt{-2}$, the conductor of $K(\theta)/K$ is 3. Assume $\theta \in L_1$, and apply Lemma 2.10 with $L = K(B[2])$, $E = K(\theta)$, $F = K$, in conjunction with the bound (2.8), to deduce that the conductor of $E/K$ is at most 2, a contradiction. For the higher layers, the 2-division tower $F_\infty$ is cyclic over $F_2$, from which it easily follows that $F_n \subset L_n \cap F_\infty \subset F_{n+1}$. Finally, consider the special case $L_2 \cap F_\infty = F_2$. Since $L_{n+1}/L_n$ has exponent 2, we can only go up one more stage at a time in the cyclotomic tower, so $L_n \cap F_\infty = F_n$ for all $n \geq 2$. ∎

## 3. Controlling the division fields

The main result of this section is the following proposition, which will be used to limit possible $\ell$-division fields for the abelian varieties under consideration.

**Proposition 3.1.** *Let $\ell$ and $p$ be distinct primes. Suppose the field $L$ satisfies the following conditions.*

(L1) $L/\mathbb{Q}$ *is a Galois extension containing* $\boldsymbol{\mu}_\ell$;
(L2) $L$ *is unramified outside $\ell$ and $p$;*
(L3) *the ramification degree of $L$ at $p$ is 1 or $\ell$;*
(L4) *if $\mathcal{D}$ denotes a decomposition group at a prime over $\ell$ in $L$, then the higher ramification groups $\mathcal{D}^u = \{1\}$ for $u > \frac{1}{\ell-1}$.*

*If $\ell = 2$ and $p = 3$ or $7$, then $L \subset \mathbb{Q}(\boldsymbol{\mu}_4, \sqrt{p})$. If $\ell = 3$ and $p = 2$ or $5$, or else if $\ell = 5$ and $p = 2$ or $3$, then $L \subset \mathbb{Q}(\boldsymbol{\mu}_\ell, p^{\frac{1}{\ell}})$.*

The first ingredient in the proof of Proposition 3.1 is a bound on the discriminant $d_{L/\mathbb{Q}}$. More generally, as Fontaine has suggested [Fo, Rem. 3.3], his methods easily imply a bound for the discriminant of the $\ell^n$-division field of a semistable abelian variety $A$ defined over $\mathbb{Q}$. See also [Jo, Thm. 2.1]. Before stating our version of the Fontaine-Joshi bound, we define the effective stage of inertia at a bad prime $p$ in this global context. Choose a prime $\mathfrak{P}$ over $p$ in $\overline{\mathbb{Q}}$ and denote by $\mathcal{D} = \mathcal{D}(\mathfrak{P}/p)$ (resp. $\mathcal{I} = \mathcal{I}(\mathfrak{P}/p)$) the decomposition group (resp. inertia group) inside $G_\mathbb{Q}$. Let $n_0$ be the minimal integer $n$ such that $\mathcal{I}$ does not act trivially on $\mathbb{Q}(A[\ell^n])$, and put $i(A, \ell, p) = n_0$. Because the inertia groups over $p$ are conjugate, $i(A, \ell, p)$ does not depend on the choice of $\mathfrak{P}$.

**Proposition 3.2.** *Suppose $A/\mathbb{Q}$ is a semistable abelian variety and let $S$ be the set of bad primes for $A$. Fix a prime $\ell \notin S$, and consider the $\ell^n$-division field $L_n = \mathbb{Q}(A[\ell^n])$. Its discriminant satisfies the inequality*

$$|d_{L_n/\mathbb{Q}}|^{\frac{1}{[L_n:\mathbb{Q}]}} < \ell^{n+\frac{1}{\ell-1}} \prod_{p \in S} p^{1-\frac{1}{\ell^n}}.$$

**Proof.** Let $\tilde{L}_n$ denote the completion of $L_n$ at a prime over $\ell$. In keeping with Fontaine's notation, we extend the valuation $v_\ell$ of $\mathbb{Q}_\ell$ to $\tilde{L}_n$ preserving $v_\ell(\ell) = 1$. By

[Fo, Thm. 1], we have

$$\text{(3.3)} \quad \frac{1}{[L_n : \mathbb{Q}]} v_\ell(d_{L_n/\mathbb{Q}}) = v_\ell(\mathfrak{D}_{\tilde{L}_n/\mathbb{Q}_\ell}) < n + \frac{1}{\ell - 1},$$

where $\mathfrak{D}_{\tilde{L}_n/\mathbb{Q}_\ell}$ is the different ideal.

Let $n_0$ be the effective stage of inertia at a prime $p \in S$. If $n < n_0$, then $L_n$ is unramified at $p$. By (2.2), for $n \geq n_0$ we have

$$\frac{1}{[L_n : \mathbb{Q}]} v_p(d_{L_n/\mathbb{Q}}) = 1 - \frac{1}{\ell^{n-n_0+1}} \leq 1 - \frac{1}{\ell^n}. \quad \blacksquare$$

Under the conditions (L1) - (L4) of Proposition 3.1, the arguments in the proof of Proposition 3.2 imply that the discriminant $d_{L/\mathbb{Q}}$ satisfies the inequality

$$\text{(3.4)} \quad |d_{L/\mathbb{Q}}|^{1/[L:\mathbb{Q}]} < \ell^{1+\frac{1}{\ell-1}} p^{1-\frac{1}{\ell}}.$$

In Table 1, we give the corresponding upper bound of Odlyzko [Od] and Diaz y Diaz [Di] on the degree $[L : \mathbb{Q}]$ for the relevant values of $\ell$ and $p$.

| $\ell$ | $\ell = 2$ | | $\ell = 3$ | | $\ell = 5$ | |
|---|---|---|---|---|---|---|
| $p$ | 3 | 7 | 2 | 5 | 2 | 3 |
| Bound for $|d_{L/\mathbb{Q}}|^{1/[L:\mathbb{Q}]}$ | 6.93 | 10.59 | 8.25 | 15.20 | 13.02 | 18.01 |
| Bound for $[L : \mathbb{Q}]$ | 10 | 22 | 14 | 68 | 40 | 168 |

TABLE 1. Bounds on the $\ell$-division field

The next ingredient in the proof of Proposition 3.1 is a class field theoretic lemma. In general, if $E \supset F \supset \mathbb{Q}$ is a tower of fields, we write $\lambda_E$ for a prime over $\ell$ in $E$ and $\lambda_F$ for $\lambda_E \cap F$. If $E/F$ is abelian, we write $\mathfrak{f}_\lambda(E/F)$ for the local conductor exponent of $E_{\lambda_E}/F_{\lambda_F}$.

**Lemma 3.5.** *Let $L$ be a field satisfying properties* (L1) - (L4), *and assume in addition that $\boldsymbol{\mu}_4 \subset L$ if $\ell = 2$. Let $F = \mathbb{Q}(\boldsymbol{\mu}_\ell)$ if $\ell$ is odd (resp. $F = \mathbb{Q}(\boldsymbol{\mu}_4)$ if $\ell = 2$). Suppose the class number of $F$ is 1, and let $s$ denote the number of primes over $p$ in $F$. Let $E$ be the maximal subfield of $L$ abelian over $F$. Then $\text{Gal}(E/F)$ is annihilated by $\ell$ and has rank at most $s$. In particular, if there is one prime over $p$ in $F$, then $E \subset F(p^{1/\ell})$.*

**Proof**. Put $A_F^\times$ for the idele group of $F$. Since $F$ has class number 1, we have

$$A_F^\times = \left(\prod_v U_v\right) F,$$

where $U_v$ is the group of units in the completion $F_v$ and $F$ is the image of $F$ on the diagonal of $A_F^\times$. Denote principal units by $U_v^{(1)}$. Write $\mathcal{N}_v$ for the image in $U_v$ of the local units of $E \otimes F_v$ under the norm map, so that the extension $E$ corresponds to the normic subgroup
$$\mathcal{N}_{E/F} = \left(\prod_v \mathcal{N}_v\right) F.$$
(By abuse of notation, write $U_v = \mathbb{C}^\times$ for the archimedean places of $F$, all of which are complex.) Let $\Upsilon = \{\nu_1, \ldots, \nu_s\}$ denote the set of primes over $p$ in $F$ and abbreviate $\lambda = \lambda_F$ for the prime over $\ell$ in $F$.

Since $E/F$ is unramified outside $\{\lambda\} \cup \Upsilon$, and the ramification degree is 1 or $\ell$ over each $\nu \in \Upsilon$, we have
$$\mathcal{N}_{E/F} \supset \left(U_\lambda^{(f_\lambda)} \prod_{\nu \in \Upsilon}(U_\nu^{(1)} U_\nu^\ell) \prod_{v \notin \{\lambda\} \cup \Upsilon} U_v\right) F.$$

It follows from the fact that the images of the (global) cyclotomic units generate $U_\lambda/U_\lambda^{(1)}$ that there is no tame ramification locally over $\lambda$ in $E/F$. In view of Fontaine's inequality (2.8) on the upper ramification numbering and Lemma 2.10, we find that $\mathfrak{f}_\lambda(E/F) \leq 2$. If $\ell$ is odd, then the image of a primitive $\ell^{th}$ root of unity generates $U_\lambda^{(1)}/U_\lambda^{(2)}$. If $\ell = 2$, then $i$ generates $U_\lambda^{(1)}/U_\lambda^{(2)}$. Thus there is a surjective map
$$\iota: \prod_{\nu \in \Upsilon} k_\nu^\times / k_\nu^{\times \ell} \to A_F^\times / \mathcal{N}_{E/F}.$$
We may conclude that $\mathrm{Gal}(E/F)$ is annihilated by $\ell$, and its rank is at most $|\Upsilon| = s$. Of course if there is no ramification over $p$ in $E$, the image of $\iota$ is trivial and $E = F$. If $p$ ramifies in $E$ and $s = 1$, there is a unique candidate for $E$, and $E = F(p^{1/\ell})$ does satisfy the desired conditions. ∎

We shall also utilize the following elementary facts from group theory.
(G1) If $|H| = 2n$ with $n$ odd, then $H$ admits a quotient of order 2.
(G2) If $H$ is an $\ell$-group and $|H| \geq \ell^2$, then $H$ admits a quotient of order $\ell^2$, necessarily abelian.
(G3) If $H$ is not an $\ell$-group and $H$ contains only one $\ell$-Sylow subgroup $S_\ell$, then $S_\ell$ is normal in $H$, and $|H/S_\ell|$ is prime to $\ell$.

Turning to the proof of Proposition 3.1, it is convenient to separate the discussion between odd $\ell$ and $\ell = 2$, although there is some overlap in the arguments.

**Proof of Proposition 3.1**, $\ell$ odd. Let $H = \mathrm{Gal}(L/\mathbb{Q}(\boldsymbol{\mu}_\ell))$. If $H$ is trivial, we are done, so we assume $H \neq \{1\}$. Let $E$ be the maximal subfield of $L$ abelian over $\mathbb{Q}(\boldsymbol{\mu}_\ell)$. For all cases in Table 1, we have $|H| < 60$, so $H$ is solvable and $\mathrm{Gal}(E/\mathbb{Q}(\boldsymbol{\mu}_\ell))$ is the (non-trivial) maximal abelian quotient of $H$. For all cases in Table 1, there is one prime over $p$ in $F$. We may conclude from Lemma 3.5, that $E = \mathbb{Q}(\boldsymbol{\mu}_\ell, p^{1/\ell})$. In particular $\ell$ divides $|H|$. But if $|H| = \ell$, then $L = E$ and we are done. In view of (G2), it now suffices to assume $H$ is not an $\ell$-group in the hope of arriving at a contradiction. By (G3) the number of $\ell$-Sylow subgroups must have the form $1 + c\ell$, with $c \geq 1$. In particular $|H| \geq \ell(1 + \ell)$.

We complete the argument with an analysis of the cases. If $\ell = 3$ and $p = 2$, then $[L : \mathbb{Q}] \leq 14$, so $|H| \leq 7$. But $|H| \geq \ell(1 + \ell) = 12$, a contradiction.

For $\ell = 5$, we may treat $p = 2$ and $p = 3$ simultaneously. At worst, we have $[L : \mathbb{Q}] \leq 168$, so $|H| \leq 42$. But we may assume $|H|$ has a divisor of the form $\ell(1 + c\ell)$ with $c \geq 1$. Then we are reduced to considering only $|H| = 30$. We rule this out by using (G1).

Suppose $\ell = 3$ and $p = 5$. Since $[L : \mathbb{Q}] \leq 68$, we have $|H| \leq 34$. Under the assumption that $|H|$ has a divisor of the form $\ell(1 + c\ell)$ with $c \geq 1$, we are reduced to considering $|H| = 12, 21, 24$ or $30$. But (G1) eliminates $|H| = 30$. For the rest of the argument, $[L : E] = 4, 7$ or $8$. Of course $L/E$ is unramified outside primes dividing 15. By property (L3), the ramification degree at each prime over 5 in $L/\mathbb{Q}$ is 3. Since $E$ already absorbs this ramification, $L/E$ also is unramified at primes over 5. Let us consider the the unique prime $\lambda_E$ over 3 in $E$. Recall that $\lambda_F$ denotes the unique prime over 3 in $F = \mathbb{Q}(\boldsymbol{\mu}_3)$. If $\lambda_E$ did not split at all in $L/E$, there would be one prime $\lambda_L$ over $\lambda_F$ in $L/F$. But then the wild ramification subgroup for $\lambda_L$ over $\lambda_F$ would be the 3-Sylow subgroup of $H$. Since the wild ramification subgroup is normal in the decomposition group, we have a contradiction of our assumption that the 3-Sylow subgroup of $H$ is not a normal subgroup.

As a consequence of this discussion, if $[L : E] = 4$ or 8, there exists at least a quadratic extension $E'/E$ which is everywhere unramified. Since the class number of $\mathbb{Q}(\sqrt[3]{5})$ is 1, we may conclude by genus theory that the class number of $E = \mathbb{Q}(\boldsymbol{\mu}_3, \sqrt[3]{5})$ is odd, in contradiction to the existence of $E'$. Suppose $[L : E] = 7$, in which case $L$ itself is an everywhere unramified extension of $E$. We may compute the absolute discriminant of $E$ as $3^7 5^4$. Since $L/E$ is unramified, the discriminant $d_{L/\mathbb{Q}}$ satisfies

$$|d_{L/\mathbb{Q}}|^{\frac{1}{[L:\mathbb{Q}]}} = 3^{7/6} 5^{2/3} \leq 10.54\,.$$

But then the Odlyzko bounds force $[L : \mathbb{Q}] \leq 22$, a contradiction. ∎

**Proof of Proposition 3.1**, $\ell = 2$. In all cases, $\mathrm{Gal}(L/\mathbb{Q})$ is small enough to be solvable. It therefore suffices to assume that the maximal subfield $E_0$ of $L$ abelian over $\mathbb{Q}$ is a non-trivial extension of $\mathbb{Q}$. By the Kronecker-Weber theorem $E_0$ is a subfield of $\mathbb{Q}(\boldsymbol{\mu}_{2^\infty}, \boldsymbol{\mu}_{p^\infty})$. But the ramification over $p$ is tame of degree at most 2, and by Proposition 2.11 we have $E_0 \cap \mathbb{Q}(\boldsymbol{\mu}_{2^\infty}) \subset \mathbb{Q}(\boldsymbol{\mu}_4)$. It follows that $E_0 \subset \mathbb{Q}(\boldsymbol{\mu}_4, \sqrt{p})$. If $L = E_0$, we are done. We may therefore assume that the maximal subfield $E_1$ of $L$ abelian over $E_0$ properly contains $E_0$. By maximality, $E_1$ is Galois over $\mathbb{Q}$.

Suppose $i \notin L$, so $E_0 = \mathbb{Q}(\sqrt{\pm p})$. Consider the maximal subfield $E_2$ of $E_1$ whose degree over $E_0$ is a power of 2. By maximality, $E_2$ is Galois over $\mathbb{Q}$. We claim that $E_2 = E_0$. Otherwise, $\mathrm{Gal}(E_2/\mathbb{Q})$ is a 2-group whose order is a least 4. It follows from (G2) that $L$ contains an abelian extension of $\mathbb{Q}$ of degree 4, contradicting the fact that $[E_0 : \mathbb{Q}] = 2$. We may conclude that if $E_0 = \mathbb{Q}(\sqrt{\pm p})$, then the degree $[E_1 : E_0]$ is odd. As a consequence of (L3), the unique prime over $p$ in $E_0$ is unramified in $E_1/E_0$. Consider any prime $\lambda_1$ over 2 in $E_1$, and let $\lambda_0 = \lambda_1 \cap E_0$ be the corresponding prime over 2 in $E_0$. Put $n$ for the ramification degree of $\lambda_1$ over $\lambda_0$, necessarily odd (tame). It follows that the residue field $k_0$ of $\lambda_0$ must contain $\boldsymbol{\mu}_n$. We now show $n = 1$, breaking up the argument according to whether $p = 7$ or $p = 3$. If $p = 7$, the residue field $k_0$ is $\mathbb{F}_2$ in all cases, so $n = 1$. If $p = 3$, the residue field $k_0$ is $\mathbb{F}_2$ unless $E_0 = \mathbb{Q}(\sqrt{-3})$, in which case we must consider $n = 3$. But class field theory or Kummer theory shows that $\mathbb{Q}(\sqrt{-3})$ does not have a Galois cubic extension unramified outside 2. At this point, we have produced a non-trivial unramified abelian extension $E_1/E_0$. But for $p = 3$ or $p = 7$, the class number of $\mathbb{Q}(\sqrt{\pm p})$ is 1, a contradiction.

It remains to study the case $i \in L$. If $E_0 = \mathbb{Q}(\boldsymbol{\mu}_4)$, then $E_1 = \mathbb{Q}(\boldsymbol{\mu}_4, \sqrt{p})$ by Lemma 3.1. This contradicts the fact that $E_0$ already is maximal abelian over $\mathbb{Q}$. Hence $E_0 = \mathbb{Q}(\boldsymbol{\mu}_4, \sqrt{p})$ and the commutator subgroup of $H = \mathrm{Gal}(L/\mathbb{Q}(\boldsymbol{\mu}_4))$ has index 2 in $H$. Since we have assumed that $L$ properly contains $E_0$, the degree $[L : \mathbb{Q}]$ is a non-trivial multiple of 4, and $H$ cannot be a 2-group by (G2). But the bound for $p = 3$ is $[L : \mathbb{Q}] \leq 10$, and we have arrived at a contradiction to complete the discussion for $p = 3$. For $p = 7$, the bound is $[L : \mathbb{Q}] \leq 22$. Using the fact that $\mathrm{Gal}(L/\mathbb{Q})$ cannot be a 2-group, we are left to consider $[L : \mathbb{Q}] = 12$ or 20. Then $L/E_0$ is a cyclic extension of degree $n = 3$ or $n = 5$. As we argued above, the extension $L/E_0$ is unramified outside 2 by property (L3). Furthermore, if a prime $\lambda_0$ over 2 in $E_0$ ramifies in $L$, then the corresponding residue field $k_0$ must contain $\boldsymbol{\mu}_n$. But $k_0 = \mathbb{F}_2$. Therefore $L/E_0$ is a non-trivial unramified abelian extension, contradicting the fact that the class number of $E_0 = \mathbb{Q}(\boldsymbol{\mu}_4, \sqrt{7})$ is 1. ∎

## 4. Non-existence results

If $A/\mathbb{Q}$ is an abelian variety with semistable bad reduction at $p$, the structure of $\mathbb{T}_\ell(A)$ as a Galois module for $G_{\mathbb{Q}_p}$ sometimes can be used to construct a $\mathbb{Q}$-isogeny that raises the effective stage of inertia $i(A, \ell, p)$ or increases the group of connected components $\Phi_{\hat{A}}(\bar{\mathbb{F}}_p)$.

**Proposition 4.1.** *Let $A/\mathbb{Q}$ be an abelian variety with semistable bad reduction at the prime $p$. Suppose $A$ has good reduction at a prime $\ell$, and assume there is one prime over $p$ in the $\ell$-division field $L = \mathbb{Q}(A[\ell])$. Then there exists a $\mathbb{Q}$-isogeny $\varphi : A \to A'$ such that $i(A', \ell, p) = i(A, \ell, p) + 1$ and $\mathrm{ord}_\ell(\Phi_{\hat{A'}}(\bar{\mathbb{F}}_p)) = \mathrm{ord}_\ell(\Phi_{\hat{A}}(\bar{\mathbb{F}}_p)) + t$, where $t$ is the toroidal dimension of the Néron fiber over $p$.*

*Suppose further that $A$ is principally polarized over $\mathbb{Q}$. If $\ell$ is odd, assume $\mathrm{Gal}(L/\mathbb{Q})$ contains an element of order 2 acting by inversion on $\boldsymbol{\mu}_\ell$. If $\ell = 2$, assume $|\mathrm{Gal}(L/\mathbb{Q})|$ is a power of 2. Then we may also arrange for $A'$ to be principally polarized over $\mathbb{Q}$.*

**Proof**. Fix a prime $\mathfrak{P}$ over $p$ in $\bar{\mathbb{Q}}$ and let $\mathcal{D}$ and $\mathcal{I}$ be the corresponding decomposition and inertia groups. Then $\mathcal{M}_1 = \mathbb{T}_\ell(A)^\mathcal{I}$ and $\mathcal{M}_2 = (\mathbb{T}_\ell(\hat{A})^\mathcal{I})^\perp$ are modules for $\mathcal{D}$.

Write $\bar{\mathcal{M}}_1$ and $\bar{\mathcal{M}}_2$ for their respective projections to $A[\ell]$. By assumption, $\mathcal{D}$ maps onto $G = \mathrm{Gal}(L/\mathbb{Q})$ by restriction. It follows that $\kappa = \bar{\mathcal{M}}_1$ or $\bar{\mathcal{M}}_2$, is a $G$-module, and therefore a $G_{\mathbb{Q}}$-module. Hence $\kappa$ is the kernel of a $\mathbb{Q}$-isogeny. We may compute the change in size of $\Phi$ using Lemma 2.4.

Under the additional assumptions, we may also arrange for $\kappa$ to be a maximal isotropic subspace of $A[\ell]$ as in Lemma 2.6 and Lemma 2.7. Then $A'$ is principally polarized.  ∎

**Theorem 4.2.** *For $p \leq 7$, there does not exist a semistable abelian variety defined over $\mathbb{Q}$ with good reduction outside the prime $p$.*

**Proof**. Suppose there exists an abelian variety $B/\mathbb{Q}$ with good reduction outside the prime $p$, and $p \leq 7$. Since Fontaine has already treated the issue of everywhere good reduction, we may assume bad reduction of semistable type at $p$. Among the finitely many abelian varieties isogenous to $B$ over $\mathbb{Q}$, we choose a variety $A$ such that the effective stage of $p$-adic inertia $i(A, \ell, p)$ is maximal. In view of Proposition 3.1, there is one prime over $p$ in the $\ell$-division field $L = \mathbb{Q}(A[\ell])$ for the values of $\ell$ and $p$ in the table of §3. Hence there exists a $\mathbb{Q}$-isogeny $A \to A'$, as constructed in Proposition 4.1. But this contradicts maximality of $i(A, \ell, p)$.  ∎

As another application when $\ell = 2$, we briefly remark on the non-existence of semistable abelian varieties $A$ over $\mathbb{Q}$, with certain "small" 2-division fields. Neumann [Ne] and Setzer [Sz] independently gave a precise description of the family of elliptic curves having good reduction outside one prime $p$ and a rational point of order 2. In particular, it is necessary and sufficient that $p$ be of the form $p = n^2 + 64$ or else $p = 17$. It seems plausible to us that if the 2-division field $L = \mathbb{Q}(A[2])$ is small, for example in the sense that $G = \mathrm{Gal}(L/\mathbb{Q})$ is a 2-group, then $A$ is isogenous to a product of Setzer-Neumann elliptic curves. We plan to address this general question in the future, perhaps with the aid of additional tools arising from the work of Schoof [Sc]. For now, we have the following limited result.

**Proposition 4.3.** *Suppose $A/\mathbb{Q}$ is a semistable abelian variety with good reduction outside one prime $p$. Let $L = \mathbb{Q}(A[2])$ be the 2-division field and assume $G = \mathrm{Gal}(L/\mathbb{Q})$ is nilpotent. Then $p \equiv 1 \pmod 4$.*

**Proof**. Quite generally, if the Galois group $G$ of the 2-division field of a semistable abelian variety is nilpotent, then $G$ is in fact a 2-group. Otherwise, there exists a maximal normal subgroup $N$ of $G$ such that $G/N$ is cyclic of odd prime order. Let $E$ be the fixed field of $N$. Since the ramification in $L$ over each odd prime of bad reduction has degree dividing 2, the extension $E/\mathbb{Q}$ must be unramified outside 2. But there is no non-trivial abelian extension of $\mathbb{Q}$ of odd degree and unramified outside 2.

We may therefore assume $G$ is a 2-group, whence $L_\infty = \mathbb{Q}(A[2^\infty])$ is a pro-2 extension of $\mathbb{Q}$. It follows that the degree over $\mathbb{Q}$ of the 2-division field of any variety isogenous to $A$ also is a power of 2. Among the finitely many varieties $\mathbb{Q}$-isogenous to $A$, choose one for which the effective stage of $p$-adic inertia acting on the 2-adic Tate module is maximal. By abuse of notation, we continue to denote this variety by $A$.

Consider the decomposition group $\mathcal{D} = \mathcal{D}(\mathfrak{p}/p)$ for a prime $\mathfrak{p}$ over $p$ in $L$. If $\mathcal{D}$ is a proper subgroup of $G$, there exists a subgroup $H$ of $G$ containing $\mathcal{D}$ such that $[G : H] = 2$. The fixed field $F$ of $H$ is a quadratic field unramified outside 2 and split completely over $p$. By the Kronecker-Weber theorem, $F \subset \mathbb{Q}(\boldsymbol{\mu}_{2^\infty})$. In fact $F \subset \mathbb{Q}(\boldsymbol{\mu}_4)$ by Proposition 2.11. Now assume $p \equiv -1 \pmod{4}$. Since $p$ is inert in $\mathbb{Q}(\boldsymbol{\mu}_4)$, we have a contradiction unless $\mathcal{D} = G$ and there is one prime over $p$ in $L$.

Assuming $p \equiv -1 \pmod{4}$, we may apply Proposition 4.1, to find a $\mathbb{Q}$-isogenous variety $A'$ such that $i(A', 2, p) = i(A, 2, p) + 1$. This contradicts maximality of $i(A, 2, p)$. ∎

## References

[Di] F. Diaz y Diaz, *Tables minorant la racine n-ième d'un corps de degré n*, Ph.D. Thesis, Publ. Math. d'Université d'Orsay, 1980.

[Ed] B. Edixhoven, *On the prime to p-part of the group of connected components of Néron models*, Compositio Math. **97** (1995), 29–49.

[Fa] G. Faltings, *Endlichkeitssätze für abelsche Varietes über Zahlkörpern*, Invent. Math. **73** (1983), 349–366.

[Fo] J.-M. Fontaine, *Il n'y a pas de variété abélienne sur $\mathbb{Z}$*, Invent Math. **81** (1985), 515–538.

[Gr] A. Grothendieck, *Modèles de Néron et monodromie*. Séminaire de Géométrie 7, Exposé IX, Lecture Notes in Mathematics **288** Springer-Verlag, New York, 1973.

[Jo] K. Joshi, *Remarks on the methods of Fontaine and Faltings*, IMRN **22** (1999), 1199–1209.

[MS] J.R. Merriman and N.P. Smart, *Curves of genus 2 with good reduction away from 2 with a rational Weierstrass point*, Math. Proc. Camb. Phi. Soc. **114** (1993), 203–214.

[Mi] J.S. Milne, *Abelian varieties*, in: Arithmetic Geometry, G. Cornell and J.H. Silverman, ed., Springer-Verlag, New York, 1986, 103–150.

[Ne] O. Neumann, *Elliptische Kurven mit vorgeschriebenem Reduktionsverhalten II*, Math. Nach. **56** (1973), 269–280.

[Od] A. Odlyzko, *Lower bounds for discriminants of number fields*, II, Tôhoku Math. J. **29** (1977), 209–216.

[Sc] R. Schoof, *Abelian varieties over number fields with everywhere good reduction*, Lecture at Durham Conference, August 2000.

[Se] J.-P. Serre, Local Fields, Lecture Notes in Math. **67**, Springer-Verlag, New York, 1979.

[ST] J.-P. Serre and J.T. Tate, *Good reduction of abelian varieties*, Ann. Math. **68** (1968), 492–517.

[Sz] C.B. Setzer, *Elliptic curves of prime conductor*, J. London Math. Soc. **10** (1975), 367–378.

[Wu] G. Wüstholz, *The finiteness theorems of Faltings*, in: Rational Points, G. Faltings and G. Wüstholz, ed., Aspects of Mathematics, **E6**, F. Vieweg & Sohn, Braunschweig, 1984.

(A. Brumer) Department of Mathematics, Fordham University, Bronx, NY 10458
*E-mail address*: brumer@murray.fordham.edu

(K. Kramer) Department of Mathematics, Queens College (CUNY), Flushing, NY 11367
*E-mail address*: kramer@forbin.qc.edu